\documentclass[12pt]{article}

\usepackage{graphicx}%
\usepackage{multirow}%
\usepackage{amsmath,amssymb,amsfonts}%
\usepackage{amsthm}%
\usepackage{mathrsfs}%
\usepackage[title]{appendix}%
\usepackage{xcolor}%
\usepackage{textcomp}%
\usepackage{manyfoot}%
\usepackage{booktabs}%
\usepackage{algorithm}%
\usepackage{algorithmicx}%
\usepackage{algpseudocode}%
\usepackage{listings}%

\usepackage{enumitem}
\usepackage{lineno}
\usepackage{float}
\usepackage{changepage}
\usepackage{aliascnt}
\usepackage{morefloats}
\usepackage{graphics}
\usepackage{graphicx}
\usepackage{listings}
\usepackage{xcolor}
\usepackage{float}
\usepackage{verbatim}
\usepackage{geometry}
\usepackage{centernot}
\usepackage[utf8]{inputenc}
\usepackage[english]{babel}
\usepackage{setspace}
\usepackage{thmtools}
\usepackage{xcolor}
\usepackage[bottom]{footmisc}
\usepackage{euscript}
\usepackage{mathtools}
\usepackage{enumitem}
\usepackage [autostyle, english = american]{csquotes}
\MakeOuterQuote{"}
\newcommand{\subscript}[2]{$#1 #2$}
\DeclareMathAlphabet{\itbf}{OML}{cmm}{b}{it}
\providecommand{\keywords}[1]{\textbf{Keywords:} #1}

\newcommand{\bea}{\begin{eqnarray*}}
\newcommand{\eea}{\end{eqnarray*}}
\newcommand{\bean}{\begin{eqnarray}}
\newcommand{\eean}{\end{eqnarray}}
\newcommand{\p}{\partial}
\newcommand{\f}{\frac}

\newcommand{\ds}{\displaystyle}
\newcommand{\no}{\nonumber}
\newcommand\ov{\overline}
\newcommand{\ri}{\rightarrow}
\newcommand{\sm}{\setminus}

\theoremstyle{thmstyleone}

\newtheorem{thm}{Theorem}
\newtheorem{lem}[thm]{Lemma}
\newtheorem{prop}[thm]{Proposition}

\theoremstyle{thmstyletwo}%

\newcommand{\doubleR}{\mathbb{R}}
\newcommand{\RR}{\mathbb{R}}
\newcommand{\CC}{\mathbb{C}}

\renewcommand{\ov}[1]{\overline{#1}}

\begin{document}

\author{Mahadevan Ganesh
\thanks{Department of Applied Mathematics and Statistics,
    Colorado School of Mines, Golden, CO,
   USA. Corresponding author email: mganesh@mines.edu},
	 Stuart C. Hawkins
	\thanks{Department of Mathematical and Physical Sciences,
    Macquarie University, Sydney, NSW 2109,
    Australia}, 
Darko  Volkov 
\thanks{Department of Mathematical Sciences,
Worcester Polytechnic Institute, Worcester, MA 01609.  
} }

 \raggedbottom

 \raggedbottom

\title{Machine learning on manifolds for inverse scattering:  Lipschitz stability analysis}
\maketitle
\abstract{%
Establishing Lipschitz stability estimates is crucial for ensuring the mathematical robustness of neural network (NN) approximations in machine learning (ML)-based parameter estimation, particularly in physics-informed settings. In this work, we derive such estimates for the inverse of a nonlinear map defined on a manifold that captures both unknown parameters and the nonlinear physical processes they influence. Our analysis is based on finite-dimensional, learnable representations of the manifold and provides 
Lipschitz stability estimates on the manifold-based subspaces,  for a class of inverse maps associated with parameter dependent linear compact operators. Such operators model scattered- and far-field data that can be used to detect structures such as cracks.

We apply our theoretical ML manifold framework to inverse Helmholtz problems in unbounded regions exterior to cracks, addressing the scattered-field data-driven inverse problem while ensuring injectivity conditions on the manifold—a requirement for the Lipschitz stability. Our method accurately recovers crack-defining parameters without requiring prior knowledge of inputs such as incident wave types or external forces on the crack. Numerical experiments using NN approximations confirm the accuracy, efficiency, and robustness of the proposed approach.}



\section{Introduction}
Machine learning (ML)-based techniques have become a mainstream computational strategy for simulating 
physical processes in bounded domains \cite{pinn_survey_2022,  thuerey2021pbdl} and also  in  unbounded regions \cite{surr_bayes_2022}. These approaches often rely on shallow or deep neural network (NN) approximations of associated learnable functions, whose existence and convergence require  regularity conditions \cite{de2021approximation, shen2021neural, yarotsky2017error}.

For example, in the ML-based solution of an ill-posed inverse problem in medical imaging studied in \cite{M-RIP-pap}, the learnability of a function defined on an image manifold requires the so-called $\mathcal{M}$-RIP condition for the governing operator.  However, while such ML models have shown strong numerical performance, the mathematical justification for Lipschitz regularity in tandem with learnable manifolds remains an open question. This article is motivated by several open research questions posed in Pages~797--798 of~\cite{alberti2023inverse}, including:
\begin{itemize}
\item Investigating unknown learnable  manifolds in conjunction with inverse problems;
\item Extending results to inverse scattering problems involving wave equations;
\item Developing numerical implementations to validate the theory.
  \end{itemize}

The motivating application in this article is the inverse wave-scattering problem of identifying a crack from scattered-field data.
In our simulations, the cracks are modeled by a parameter $m$ which lies in a finite-dimensional space.
The scattered field $u$ is governed by the Helmholtz partial differential equation (PDE) in an unbounded region exterior to the crack. This field may in practice
result from an \textbf{unknown} forcing term.
This forcing term may be 
associated with an unknown incident wave
interacting with the crack; it could also be any kind of unknown source 
$\varphi$
for $u$ supported on the crack.

Next, we briefly describe our general framework based  main theoretical and computational results in this article: Let $A_m$ be the linear operator mapping $\varphi$ to the restriction of $u$ to $S_R$ where, for
example, $S_R$ is a circle
 with radius $R$ encompassing the region where the crack may lie. 
 Our main result is Theorem \ref{main theorem}, which states the inverse Lipschitz stability estimate for a class of operators:
\bea
\| A_m \varphi - A_{m'}  \psi  \| \geq 
C(\| \psi\| |m-m'| + \| \varphi -\psi\|),
\eea
where $m, m'$ are two geometry parameters and $\varphi, \psi$ are two forcing terms that are  linear combinations 
of the first $N$ singular functions of $A_m$ and $A_{m'}$.
Theorem  \ref{main theorem} (proved in the Appendix) is actually derived in a general theoretical framework that
can be applied to 
applications beyond the one presented in this paper.
The practical implication of Theorem
\ref{main theorem} is that the inverse function  $ A_m \varphi \mapsto (m , \varphi) $ 
has a Lipschitz continuous inverse on the set
\bea
 \{ A_m  \varphi: m \in {\cal B},  \varphi \in E_{m,N}  , \| A_m \varphi \| =1 \}, 
\eea
where ${\cal B}$ is the set of possible parameters $m$ and 
$E_{m,N}$ is the space of linear combinations of the first $N$ singular functions of $A_m.$
It is then known that this inverse can be 
approximated by an NN. The depth of the NN and of the number nodes
required to achieve given accuracy can be estimated~\cite{shen2021neural, yarotsky2017error,de2021approximation}.

In the inverse scattering application presented in this paper, $A_m \varphi$ corresponds to measurements of the  
scattered  field taken on a curve surrounding the scatterer. The fact  that the forcing term $\varphi $ is also unknown allows us to
model passive inverse problems.
Two conditions,  denoted \ref{U1}--\ref{U2}, are required for Theorem  \ref{main theorem}  to hold.
They correspond to injectivity conditions for the map $(m , \varphi) \mapsto A_m \varphi$ 
and its differential, as explained in Proposition 
 \ref{U1 U2 eq} and its proof in the Appendix. 
 For the scattering application presented in this paper, we state 
 in Proposition  \ref{BU2 cond}
 and prove in the Appendix that conditions \ref{U1}--\ref{U2}
 are satisfied. In section \ref{num section}, we describe how we compute an NN for this scattering application.
 The NN training set involves random choices for $m$ and random linear combinations of the first five singular vectors of $A_m$.
 We show 
 how this NN performs on many realizations of the data, possibly noisy.
 In this performance test, if the forcing term is due to an incoming wave, $\varphi$
  is related to the incoming wave through an integral equation. 
  In
 section    \ref{finite el section}, we make sure not to commit the so called
 inverse crime by computing data using a
completely independent numerical solver for a
slightly different model problem in which the thin crack
is modeled by a rectangle. This new model holds in a bounded domain.
This data is then fed into our NN and $m$ is successfully recovered.

 One can find in the literature many  studies on Lipschitz  stability for inverse problems,
 as discussed in the extensive survey \cite{alberti2023inverse}.  It is well known that PDE-based inverse problems are exponentially ill-posed  \cite{alessandrini1997examples, alessandrini2009stability, hadamard2014lectures}. Nonetheless, under restrictive assumptions—most often involving finite-dimensional subspaces—Lipschitz stability can be  established.
However, these works typically do not address NN-relevant learnability conditions such as those in \cite{M-RIP-pap}. 

For instance, the inverse conductivity problem has been shown to be Lipschitz stable when restricted to a finite-dimensional space of piecewise constant conductivities \cite{de2012local, harrach2023calderon}. In some other simplified models, the medium is assumed to consist of a known background and an unknown inclusion with different conductivity. Under the assumption that the inclusion is polygonal, global Lipschitz stability has been proved \cite{beretta2022global}. Similarly, in \cite{triki2019stability, volkov2021stability}, inverse problems for recovering cracks in elastic or Laplacian-governed media were shown to be Lipschitz stable, provided the crack geometry can be parameterized by a finite-dimensional vector.

We also recall that many studies derive Lipschitz regularity estimates in terms of Dirichlet-to-Neumann (DtN)\
 boundary operators \cite{alessandrini2005lipschitz, beretta2013lipschitz, beretta2022global, harrach2019uniqueness, harrach2023calderon}. However, such estimates are often impractical, as they assume access to the full DtN operator, which effectively requires the ability to impose infinitely many boundary conditions. In contrast, our work focuses on passive inverse problems, where only a single, and typically unknown, forcing term is available.

Some recent work has studied inverse Lipschitz regularity for general nonlinear inverse mappings between Banach spaces
  \cite{alberti2023inverse, alberti2022infinite}. A key idea in these studies is to restrict the nonlinear map to finite-dimensional subspaces. If both the restricted map and its differential are injective, then by the inverse function theorem, one can establish a Lipschitz homeomorphism.  In our work,
we use Theorem~2.2 from \cite{alberti2023inverse} to derive our main theoretical result. This greatly simplifies the proof 
of that result. Theorem~2.2 from \cite{alberti2023inverse} provides a local version of our main result where $m$ is reduced 
to a small neighborhood. In the proofs provided in the Appendix, we explain how to obtain a result which is global in $m$.

\section{Lipschitz regularity results: Manifolds and  subspaces
 }\label{prelim}
\subsection{Abstract framework and injective equivalent assumptions}\label{abs:framework}
For clarity, all the proofs of the theoretical statements from this section  are provided in the appendix section. \\
Let 
\begin{itemize}
\item $E$ and $F$ be two Hilbert spaces;
\item $A_m: E \to F$ a compact linear operator depending on a vector parameter $m$;
\item the vector parameter $m \in {\cal B}'$, an open subset of $\RR^p$;
\item 
 ${\cal B} $ a compact subset of  ${\cal B}'$;
\item the function
	$m \mapsto A_m$ be of class $C^1$ in ${\cal B}'$.
\end{itemize}

It will be useful to consider the function
\bean \label{PSsi def}
{\bf \Psi}:  {\cal B}' \times E\to F, \no \\
{\bf \Psi}(m, \varphi) = A_m \varphi.
\eean
Clearly  ${\bf \Psi} \in C^1(  {\cal B}' \times E, F)$ due to our assumptions on
$A_m$.
Next, we wish to apply an adequate version of the implicit function theorem to ${\bf \Psi}$.
Due to the particular  definition of ${\bf \Psi}$, this function is linear in its second argument. 
As a result, there is a convenient equivalent way of requiring ${\bf \Psi}$ and its differential 
to be injective. This is the object of the following Proposition.
\begin{prop}\label{U1 U2 eq}
Let ${\bf \Psi}$ be defined by \eqref{PSsi def}.
${\bf \Psi}$ is injective on ${\cal B}' \times (E \setminus \{0\})$ and  
the derivative of ${\bf \Psi}$ is injective at every point in ${\cal B}' \times (E \setminus \{0\})$
 if and only if the following two conditions hold:
\begin{enumerate}[label=(\subscript{\mathcal{E}}{{\arabic*}})]
\item
\label{U1} For any $m,m' \in {\cal B}'$, 
for any $\varphi, \psi \in E$, if $\varphi\neq 0$ and $A_m  \varphi = A_{m'}\psi$
then $m=m'$ and $\varphi=\psi$.
\item \label{U2} 
For $q \in \RR^p$ with $|q|=1$, denote by $\p_q A_m$ the derivative of $A_m$
in $m$ in the direction of $q$. 
The linear operator $T: E \times E \ri F$, given by $T(\varphi, \psi) = \p_q A_m \varphi+ A_m \psi$,
is injective.
\end{enumerate}
\end{prop}

Before stating Lipschitz stability  results for the  inverse of ${\bf \Psi}$
(under adequate assumptions), we first
provide simple examples
illustrating how  assumptions
\ref{U1}-\ref{U2} can hold in practice.
Later, in section~\ref{inv_crack_example},  we will verify assumptions 
\ref{U1}-\ref{U2}
for a class of
wave propagation models. 

\vspace{0.1in}
{\bf Generic Example-1 (Finite Dimensional spaces $E$ and $F$)}:

Choose ${\cal B} = [1,2]^2$, $E=\RR^n$ with its natural basis
	$e_1, ..., e_n$,
	$F=\RR^{3n}$ with its natural basis
	$f_1, ..., f_{3n}$.
	For $m=(m_1, m_2) \in {\cal B} $ for our first generic example, 
	define $A_m$ by setting for $j=1, ..., n$,
	$$A_m e_j = m_1 f_j + m_2 f_{j+n} + m_2^2 f_{j + 2n}.$$
	
	To verify \ref{U1} for this example, let  $\varphi=\sum_{j=1}^n \varphi_j e_j $, $\psi=\sum_{j=1}^n \psi_j e_j \in E$
	such that $\varphi \neq 0$. 
	Let $m, m' \in {\cal B}$. Assume that $A_m \varphi = A_{m'}\psi$.
	Then for some $k \in \{ 1, ..., n\}$, $\varphi_k \neq 0$
	and
	$m_1 \varphi_k = m_1' \psi_k$, $m_2 \varphi_k = m_2' \psi_k$, 
	and $m_2^2\varphi_k = m_2'^2 \psi_k$.
	This implies that $m_2= m_2'$, so $\varphi_k=\psi_k$ and $m_1=m_1'$.
        The result
	$\varphi=\psi$ easily follows.
	
Next to verify the assumption \ref{U2},	we let $q=(\cos \alpha, \sin \alpha)$
and assume that for $\varphi=\sum_{j=1}^n \varphi_j e_j $, $\psi=\sum_{j=1}^n \psi_j e_j \in E$
and some $m \in {\cal B}$,
	$\p_q A_m \varphi + A_m \psi=0$. 
	Then for $j=1, ..., n$,
	\bea
	(\cos \alpha) \varphi_j + m_1 \psi_j = 0, \\
	(\sin \alpha) \varphi_j + m_2 \psi_j =0, \\
	2 m_2 (\sin \alpha) \varphi_j + m_2^2 \psi_j =0.
	\eea
	If $\sin \alpha =  0$, then the first  two equations imply $\varphi_j = \psi_j =0$.
	If $\sin \alpha \neq 0$, then the last two equations imply $\varphi_j = \psi_j =0$.
	In all cases, we found that $\p_q A_m + A_m$ is injective on
	$E \times E$.
	
\vspace{0.1in}
{\bf Generic Example-2 (Infinite Dimensional spaces $E$ and $F$)}:

For this example,  we choose again  ${\cal B} = [1,2]^2$.
Let $E$ be a Hilbert space and $\{ e_n: n \geq 1 \}$ a  Hilbert basis of $E$.
Let $F$ be another Hilbert space and $\{ f_n: n \geq 1 \}$ a  Hilbert basis of $F$.
	For $m=(m_1, m_2)$ in ${\cal B} $ we 
	define $A_m$ by setting for $n \geq 1$
	$$A_m e_n = \f{m_1}{n} f_{3 n }+ \f{m_2}{n} f_{3n+1} +
	\f{m_2^2}{n} f_{ 3n+2}.$$
	Note that this definition ensures that $A_m $ is compact. Similar to  calculations in  Example-1,
	it can be shown that the assumptions \ref{U1}-\ref{U2} also hold for this example.\\
	Our main theoretical result 
is a  stability result  involving  on the first singular vectors of
$A_m$. This is important in practice since  
first singular vectors can be efficiently computed  using SVD algorithms.
We now define  more precisely what is meant by  first singular vectors.
Fix a positive integer $N$. For $m \in {\cal B}$, define a subspace $E_{m,N}$
of $E$, with $\dim E_{m,N} =N $,  such that 
\bean \label{EmNdef}
 E_{m,N} 
\subset \mbox{Ker } (A_m^* A_m - \lambda_1^2(m) I ) , \no \\
\mbox{or for some integer $r$, } \no  \\
\sum_{j=1}^{r-1} \mbox{Ker } (A_m^* A_m - \lambda_j^2(m) I ) 
\subset E_{m,N} 
\subset 
\sum_{j=1}^{r} \mbox{Ker } (A_m^* A_m - \lambda_j^2(m) I ) .
\eean 
With this definition of $E_{m,N}$ we can now state our main theoretical result. 

\begin{thm} \label{main theorem}
Assume that  $A_m$ 
satisfies \ref{U1}--\ref{U2}. 
Fix a positive integer $N$ and consider the subspaces $E_{m,N}$ defined in \eqref{EmNdef}. 
Then there is a positive constant $C$ such that for all
$m, m' \in {\cal B}$ and all $\varphi \in E_{m,N}$ and $\psi \in E_{m',N} $,
\bean \label{stability EmN}
\| A_m \varphi - A_{m'}  \psi  \| \geq 
C(\| \psi\| |m-m'| + \| \varphi -\psi\|).
\eean
\end{thm}

Let us now explain the practical importance of Theorem
\ref{main theorem}.
Under the assumptions of this theorem, ${\bf \Psi}$ defined by \eqref{PSsi def}
has a Lipschitz continuous inverse on the set
\bean
 {\bf \Psi} \left( \{  (m, \varphi): m \in {\cal B},  \varphi \in E_{m,N}  , \| A_m \varphi \| =1 \}\right), \label{red_set}
\eean
  since for some constant $C>0$, $\| A_m \varphi \| \geq C \| \varphi \|$ for all
  $\varphi \in E_{m, N}$ and $m \in {\cal B}$.
It is then known that this inverse can be 
approximated by an NN.
Moreover, the depth of the NN and of the number nodes
required to achieve given accuracy can be estimated.
Indeed, there are  many papers in the NN literature
that provide upper bounds for the size of neural networks approximating
 Lipschitz functions. For example, we  refer to  
 \cite{shen2021neural, yarotsky2017error} for estimates valid if
the ReLU  function is used for activation and 
\cite{de2021approximation} if the hyperbolic tangent function ($\tanh$)  is used instead.

\section{Application: PDE and  parametric inverse functions}\label{Dirichlet screen}

In this section, we first state  the well-posedness of a specific class of forward scattered-field models governed by the Helmholtz PDE in unbounded regions of $\mathbb{R}^d$, where $d = 2, 3$. 
We then focus on a particular subset of these models in two dimensions and describe the parameter vector $m$, which characterizes the geometry of a line crack. 
That way, we define a specific   operator $A_m$ and a function ${\bf \Psi }$
such that  Theorem
\ref{main theorem} can be applied. 
There only remains to verify that this specific $A_m$ satisfies conditions \ref{U1}--\ref{U2}.
Here too, for clarity, all proofs are provided in the appendix.

\subsection{Helmholtz forward models in unbounded regions}
Let $\Gamma$ be a Lipschitz  open surface in $\RR^{d}$ if $d=3$, or
a Lipschitz open curve if $d=2$.
  Let
	$D$ be a  domain in $\RR^{d}$ with boundary $\p D$ such that
	  $\Gamma \subset \p D$.
	The trace
	theorem (which is also valid in Lipschitz domains~\cite{ding1996proof},
	allows us to define an inner and outer trace in $H^{\f12}(\p D)$ 
	of functions defined in $\RR^{d} \setminus \p D$ with local $H^1$ regularity.
	We assume that the spatially-dependent wavenumber $k$ is in $L^\infty(\doubleR^{d})$ and satisfies
	\begin{itemize}
\item $k$ is real-valued; 
\item there is a positive constant $k_{min}$ such that $k \geq k_{min}$
almost everywhere in $\doubleR^{d}$; 
\item there exists positive constants  $R_0, k_0$ such that if $|x|\geq R_0$, and $x \in 
\doubleR^{d}\sm \ov{\Gamma}$,
$k(x) = k_0$. 
\end{itemize}

        We impose a Dirichlet condition on the total field on $\Gamma$.
       The  resulting  data 
 could then be  derived from 
an incoming incident wave while the problem is solved for a scattered wave.
in which case $\Gamma$ is often called a 
screen. One can find many references in the literature for the case where 
$k^2$ is constant in space~\cite{ hsiao1991dirichlet, stephan1985augmented, wendland1990hypersingular}. In particular, these references
include an analysis of singularities of the solution at the tip of the crack and the 
analysis of numerical methods for solving these problems using integral equations
on $\Gamma$. 

Following Section~2.6 of \cite{nedelec2001acoustic}, with $r(x) = |x|,~x \in \RR^{d}$, we consider the solution space
 ${\cal V}$, defined first for $d=3$ as 
  \bea
{\cal V} = \left\{ v \in H^1_{loc}(\RR^{3}\setminus \ov{\Gamma}): \f{v}{\sqrt{1+r^{2}}}, 
\f{\nabla v}{\sqrt{1+r^{2}}}, \f{\p v}{\p r} - i k_0 v \in L^2(\RR^{3}\setminus \ov{\Gamma}) \right\}.
\eea
Its counterpart for the $d=2$ case is:
\bea
 {\cal V} =   \left\{ v \in H^1_{loc}(\RR^{2}\setminus \ov{\Gamma}): \
 \f{v}{\sqrt{1+r}
\ln(2+r)}, 
\f{\nabla v}{\sqrt{1+r} \ln(2+r)}, \f{\p v}{\p r} - i k_0 v \in L^2(\RR^{2}\setminus \ov{\Gamma})
\right\}.
\eea

We consider the following problem:  
  find $u \in {\cal V}$  such that 
	\bean
        (\Delta + k^2 )u&=&0\text{ in }\doubleR^{d}\sm \ov{\Gamma},  \label{D1}     \\
      u 
			&=&g \mbox{ on }\Gamma,  \label{D2}\\
 \f{\p u}{\p r} - i k_0 u &=& O\left(|x|^{-(d+1)/2}\right), \quad \text{as} \quad r \ri \infty
     \label{D3},
\eean
where $g$ is the restriction to $\Gamma$ of a function in 
$H^{\f12} (\p D)$,
and ${\cal V}$ is a function space on which the BVP with the radiation condition~\eqref{D3} is well posed, and whose
solution $u$ depends continuously on $g$.
In applications,  $g$ is often the negative of the trace on $\Gamma$ of the
incoming
incident field and  $u$ represents the scattered field, while
the trace of the total field vanishes on $\Gamma$.  
The following proposition  states that the crack forward  problem is uniquely solvable, which in our framework
makes it possible to later define an operator $A_m$.
\begin{prop} \label{direct  dir prob}
The BVP (\ref{D1})--(\ref{D3}) is uniquely solvable.
The solution 
$u \in {\cal V}$  depends continuously on the forcing term
$g \in H^{\f12} (\Gamma)$.
\end{prop}


If $k^2$ is constant, 
the solution $u$ to the BVP (\ref{D1})--(\ref{D3}) can be written in integral
form. 
\bean \label{Dir int}
u(x) = \int_\Gamma \Phi(x,y) \left[ \f{\p u}{\p n} (y)\right] d \sigma(y),
\eean
where $\Phi $ denotes the free space Green function for the Helmholtz equation
and  $\left[ \f{\p u}{\p n} \right] $ denotes
the jump of $\f{\p u}{\p n} \ $ across $\Gamma$,
which  is in $H^{-\f12}(\Gamma)$,
see \cite{stephan1984augmented}. 
Referring to the BVP (\ref{D1})--(\ref{D3}), suppose that
$\left[ \f{\p u}{\p n} \right] $ is zero on $\Gamma \cap B$, where $B$ 
is an open ball with center on $\Gamma$. Then $u$ is locally $H^1$
in $B$ and  satisfies $ (\Delta + k^2 )u=0$ in $B$. Then $B$ can be
taken out from $\Gamma$ without changing the solution $u$. 
We thus make the following minimal assumption on $\Gamma$: \\
\begin{center}
$\ds \left[ \f{\p u}{\p n} \right] ${\sl has full support in }$\Gamma$,
	\end{center}
or equivalently,
for any open ball $B$ centered on $\Gamma$, $u$ cannot be extended to 
	a function satisfying  $ (\Delta + k^2 )u=0$ in $B$.\\
	The following theorem states that the crack inverse problem is uniquely solvable, which in our framework
is related to condition \ref{U1}.
\begin{thm}
\label{InverseProblemResultDir}
For $i=1,2$,
    let $\Gamma_i$   be a Lipschitz open surface, let 
 $u^i$ be the unique solution to the BVP  (\ref{D1})--(\ref{D3}) with $\Gamma_i$ in place of $\Gamma$ and the Dirichlet condition $g^i \in H^{1/2}(\Gamma_i)$ in place of $g$.  
    Let $R \geq R_0$ and 
	 $S_R$ be a sphere of radius $R$. 
	Assume that  $\RR^{d} \setminus \ov{\Gamma_1 \cup \Gamma_2}$ is
	connected and that $g^i$ has full support in  $\ov{\Gamma_i}$, $i=1,2$. 
	If $u^1=u^2$ on $S_R$, then 
	$\ov{\Gamma_1}=\ov{\Gamma_2}$ and $g^1=g^2$ almost everywhere.	
\end{thm}

\subsection{Crack inverse problem: 
verifying conditions \ref{U1}-\ref{U2}}\label{inv_crack_example}
We start from the solution $u$ to the BVP  (\ref{D1})--(\ref{D3})  written in integral
form \eqref{Dir int} for the case $d=2$  with a concrete parametric description of a linear crack $\Gamma$.
For this case, the kernel in~\eqref{Dir int} is given by 
\bean \label{Phi def}
\Phi(x,y) = \f{i}{4} {\cal H}^1_0(k|x-y|).
\eean

The line supporting the linear crack $\Gamma$ can be parametrized
by choosing a unit vector direction $\tau$ and an offset scalar
parameter $a$ such that this line goes through the point
$a n \in \RR^2$, with  $n$ a unit vector normal  to $\Gamma$.
 The solution $u$ to problem (\ref{D1})--(\ref{D3}) can then be written in integral
form as
\bean
u(x) = \int_{-M}^{M} \Phi(x,y(t)) \left[ \f{\p u}{\p n} \right] (y(t)) dt, \label{bounded dir}\\
y(t) =\tau t +a n, \label{bounded2 dir}\\
\tau = (\cos \theta, \sin \theta), \quad
n=(- \sin \theta , \cos \theta),\label{bounded3 dir}
\eean 
where  $M$ is  such that the support of $\left[ \f{\p u}{\p n} \right](y)$ is in 
$[-M, M]$ with $y=\tau t +a n$.

Accordingly, with vector parameter $m = (\theta,a)$, we define the
forward model operator $A_m$ discussed in the first three sections of this article, using the framework in subsection~\ref{abs:framework} with 
$E= H^{-\f12 } ((-M,M))$,~$F= L^2(S_R)$, and for a small $\epsilon >0$,
\bean
{\cal B}= \{ (\theta, a ): -\pi/2 \leq \theta \leq  0,
 -a_{\mbox{max}} \leq  a \leq   a_{\mbox{max}} \} , \label{Mrange} \\
{\cal B}'= \{ (\theta, a ): -\pi/2 - \epsilon < \theta  < \epsilon ,
 -a_{\mbox{max}} - \epsilon < a <  a_{\mbox{max}} + \epsilon 
 \}  , \no 
\eean
where
\bean
\mbox{
the constants
$R, M, a_{\mbox{max}}
$ are such that the distance from 
the line segment}  &\no \\
 \mbox{$t \tau  + a n $, $-M \leq t \leq M$,
to the circle $S_R$ is bounded by a positive constant.} &\label{B1} 
\eean
We now focus on the inversion of the operator 
\bean \label{particular psi}
{\rm \bf \Psi} : {\cal B} \times H^{-\f12 } ((-M,M)) \ri L^2(S_R),
\eean
 induced by the solution in~\eqref{bounded dir}:
\bean
A_m \psi= u|_{S_R}, 
\mbox{ where } u(x) = \int_{-M}^{M} \Phi(x,y(t)) \psi(t) dt.\label{Bmhere2} 
\eean
It is clear due to formulation \eqref{bounded dir}
that $A_m$ is $C^1$ in $m=(\theta,a)$.
We now want to apply Theorems \ref{Alberti thm} and 
\ref{main theorem} to ${\rm \bf \Psi}$ defined in 
\eqref{particular psi}.
Theorem \ref{Alberti thm} 
was formulated for real
Banach spaces.
We now consider modifications required to accommodate
$E =H^{-\f12 } ((-M,M)) $, $F= L^2(S_R)$ being complex Hilbert spaces.
 Observe that $E,F$ and the finite dimensional spaces $R(P_{m,i})$
 defined in the proof of Proposition \ref{proj case} are also vector spaces over 
 $\RR$.  Accordingly,  Theorem \ref{Alberti thm} can be applied, 
 estimate \eqref{Wmi} holds, and the rest of the proof of Proposition \ref{proj case} can be carried out to obtain estimate
 \eqref{stability proj case} and then estimate \eqref{stability EmN}.\\
At this stage, it suffices to  prove that conditions \ref{U1}--\ref{U2} are satisfied by the  operator $A_m $ defined in \eqref{Bmhere2},
to then claim the Lipschitz regularity for the inverse of the adequately restricted operator.


\begin{prop} \label{BU2 cond}
{{\bf [Conditions \ref{U1} and \ref{U2} hold]}} \\
For all $\psi, \phi \in  H^{-\f12 }((-M,M))$, for all $m, m' \in {\cal B}'$,
if $\psi \neq 0$,
$
A_m \psi = A_{m'} \phi$  implies $m=m'$ and
$\psi =  \phi$.\\
Let  $q \in \RR^2$ be a unit vector. 
If 
$\p_{q} A_m \psi =A_m \phi $ then $\psi=\phi=0$.
\end{prop}

We can now apply Theorem \ref{main theorem} for any integer
 $N$ 
involved in the definition of the subspaces $E_{m ,N}$.
However,  as $N$ grows large the constant $C $ in 
 \eqref{stability EmN} tends to zero. 
In fact, $C=O(\tau^{-N})$ for some $\tau \in (0,1)$. 
This is due to the fact that $\Phi(x,y)$ is analytic in $y$ if
$y$ is in some open neighborhood of all possible line
segments   $t \tau  + a n $ such that condition \eqref{B1} holds.
 The exponential 
decay   $C=O(\tau^{-N})$ was proved in
 \cite{birman1977estimates, volkov2024optimal}.

\section{Numerical simulations}\label{num section}

 We now show  simulations in relation to the    inverse problem 
introduced  in
 section~\ref{inv_crack_example}.
The goal is to recover the  parameter $m \in {\cal B}$ 
from the data $A_m \psi$ defined in \eqref{Bmhere2} using an NN.
The simulations comprise three stages. 
First, data for training the NN is created and stored. 
 Second, the NN is trained on that data.
Third, the accuracy and the computational speed of the NN is tested
on entirely new data produced by incoming waves, point sources, or
arbitrary forcing terms altogether.

\subsection{Specific values of parameters and bounds}
To establish the parameter set  ${\cal B}$  defined in  \eqref{Mrange}  we fix 
$a_{\mbox{max}} =1$. In our simulations, the constant wavenumber
$k$ is 1.5 and $R$, the radius of $S_R$,  is 4.
Next we set bounds for the support of $\psi$
for condition \eqref{B1} to hold. Recalling the $y$ dependency on
$t$ in~\eqref{bounded2 dir},  we require   the support of $\psi$
with regard to $t$ in 
\eqref{Bmhere2} to be such that $t$ is in the interval with center
$o$ in [-1,1] and length $l$ in $[1, 3]$.  With these numbers the distance 
from the support of $\psi$ to $S_R$ is bounded below by $\sim 1.3$. 
 

\subsection{Discrete approximation of $A_m $ and learning data setup} \label{Discrete approximation}
In recent years,  solvability of bounded-domain PDEs have been  explored using NN-based approximations~\cite{pinn_survey_2022} including for  inverse model counterparts~\cite{belkouchi2025learning, le2024unfolded} that are often aided by automatic differentiation capabilities in frameworks such as PyTorch and TensorFlow.   These PDE informed NN  (PINN) frameworks typically incorporate the PDE directly into the NN loss function by sampling within the bounded domain, which is not practical for  our unbounded region model and associated boundary integral solution representation. 

The survey~\cite{pinn_survey_2022} emphasizes (see the last line of
its abstract) that  NN theoretical challenges remain unresolved even within the context of bounded-domain PDEs. In addition, in these 
bounded PDE based NN approaches~\cite{pinn_survey_2022}  a new NN is computed for each new instance of the boundary conditions  
since these boundary conditions are involved in the loss function. Solving inverse imaging problems through minimization can be done
using proximal neural networks (PNNs), see \cite{le2024unfolded}.
There too, a new neural network is computed for each new imaging problem.

Our computational approach is strikingly different for the unbounded domain PDE model (with established mathematically robust analysis for NN approximations). We first fix a wavenumber $k$,  a likely 
range for the crack geometry parameters $a, \theta, o, l$, a circle $S_R$, 
 and $N_S$ points on $S_R$ where the scattered field $u$ will be measured.
Denote ${\bf u}$ for associated field measurement, a vector in $\CC^{N_S}$.  
The  neural network that we compute approximates the map ${\bf u} \mapsto ( a , \theta)$.
In particular, we do not determine  a new NN for a new measurement
 $\bf{u} $. 
 The same NN can be used to recover  any crack within the likely range, for any (unknown) Dirichlet
 condition on that crack, as we shall demonstrate. 
 
Recall that 
$
A_m \psi =  \int_{-M}^{M} \Phi(x,y) \psi (t) dt$
where $x \in S_R$ and  that $t$ and $y$ 
are related by \eqref{bounded2 dir}.
 As $x$ remains bounded from $y$, $\Phi(x,y)$ is smooth.
We  then approximate the smooth function 
 $A_m \psi $ on $S_R$
by the vector $(A_m \psi (x_i))_{1 \leq i \leq N_S}$,  with
observation points $x_i = \left(R \cos (i\f{2 \pi }{N_S}), 
R\sin (i \f{2 \pi }{N_S})\right)$ evenly spaced on the circle of radius R.

Given the support of $\psi$, we have
$ \int_{-M}^{M} \Phi(x,y) \psi (t) dt =
\int_{M_1}^{M_2} \Phi(x,y) \psi (t) dt$
where $M_1 = o - l/2$, $M_2 = o + l/2$.
We then use the change of variables
$t=\f{M_2-M_1}{2} s +\f{M_2+M_1}{ 2}$, $ -1 \leq s \leq 1$ to
evaluate the
integral as $ \int_{-1}^{1} \Phi(x,y) \psi (t) \f{M_2-M_1}{2} ds$.
If $\psi$ solves the equation
\bea
    \int_{-1}^1  \Phi (x, y(t)) \psi (t) \f{M_2-M_1}{2} ds= - u_{inc}(x), \quad x \in \Gamma,
\eea
where $u_{inc}$ is an incoming wave,
it is well-known that $\psi$
presents a square root singularity at each endpoint 
of its  support \cite{stephan1984augmented}.
This
allows us to write
$\psi (t) = \tilde{\psi}(s)/\sqrt{1 -s^2}$,
where $\tilde{\psi}$ is smooth.
This motivates the change of variables $s= \sin v$. 
The integral is then approximated 
 by a finite sum
using $N_\Gamma= 10 $ quadrature points for $v$ forming
a uniform grid of $[- \pi/2, \pi/2]$.
We denote by $y(v_j)$
the associated values for $y$, with
$v_j = -\f{\pi}{2} + (j-1)\f{\pi}{N_\Gamma -1}$ for
$j=1, ..., N_\Gamma$. 
Note that $y(v_j)$ depends on $o, l$, and $m$.

Altogether, $A_m$ is approximated  by an $ N_S \times  N_\Gamma$
complex matrix 
$A_{m,app}$ with entries
\bea
\tau_j \Phi (x_i, y(v_j)) ) \f{M_2-M_1}{2} \f{\pi}{N_\Gamma -1} , \, 
i=1, ..., N_S, \,  j=1, ..., N_\Gamma.
\eea
Here the weights 
$\tau_1=\tau_{N_\Gamma}=\f12$ and $\tau_j=1, \mbox{ for } 1 < j < N_\Gamma$,
come from the trapezoidal rule, and
the constant term $\f{M_2-M_1}{2} \f{\pi}{N_\Gamma -1}$
can be omitted since we are only interested in computing singular vectors.
The bridge between application of Theorem~\ref{main theorem}
to $A_{m,app}$ instead of $A_m$ is
covered in Theorem~4.2 of~\cite{volkov2024stability}, which
asserts that estimate \eqref{stability EmN}
applies to $A_{m,app}$  as well with $C/2$ in place of $C$, 
as long as the dimension $N_S$ is sufficiently large.
In fact, this theorem implies that any numerical method based on convergent quadratures
could be used for approximating  $A_{m,app}$ from $A_m$.

We then follow these steps to produce learning data:
\begin{enumerate}
\item A random geometry is chosen for $\Gamma$. This is done 
by sampling values for $a, \theta$ in ${\cal B} $ at random using a uniform probability distribution.
\item A random support is chosen for $\psi$. This is done 
by picking random values for $o, l$ using uniform probability distributions.
\item For these choices, the matrix $A_{m,app}$ is constructed as described above.
\item 
A reduced singular value decomposition of $A_{m,app}$  is computed.
Let $v_1, ... , v_5$ be the corresponding first five eigenvectors of $A_{m,app} A_{m,app}^*$
in $\CC^{N_S}$.
\item A vector $r_1, ..., r_5$ is chosen at random in the ball of $\CC^5$ 
centered at the origin and with radius 1.
\item Let $w= r_1 v_1 + ... + r_5 v_5$.
The input for learning is  formed from separating the real and  the imaginary
parts the normalized vector
$w / \| w \|$. It is therefore a vector in $\RR^{2 N_S}$.
 The target is
$(a,\theta)$.
\end{enumerate}
The set of learning data in our NN simulations comprised $10^6$ input-target pairs.

\subsection{Neural network training step}\label{network training}
In our numerical solutions, we used the value $N_S =40$.
Accordingly, our input for learning is in $\RR^{80}$.
The targets $(\sin \theta, a)$ are in $\RR^2$.
Inputs and targets were  generated as described in the previous paragraph.
We trained a neural network ${\cal N}_1$ 
composed of 
an entry layer with  width 80, three hidden layers with width 80,
 and one exit layer with  width 2 on this data set. 
The activation function connecting these layers was chosen to be
the hyperbolic tangent function.
The training was performed using the ADAM algorithm 
\cite{kingma2014adam}.
We found that this stochastic minimization algorithm
is particularly efficient given 
the size of our problem: 
this algorithm made it possible to 
compute gradients of the penalty function 
on randomized mini-batches.  
The learning rate was set to 0.001, the mini-batch size was $10^4$ 
and the $l^2$ regularization constant was set to $10^{-5}$. More precisely,
if the weights are ${\cal W}_j$ the regularization term is $10^{-5} \ds \sum_{j}{\cal W}_j^2$.

\begin{figure}[H]
  \centering
      \hspace{-0.4in} \includegraphics[scale=.49]{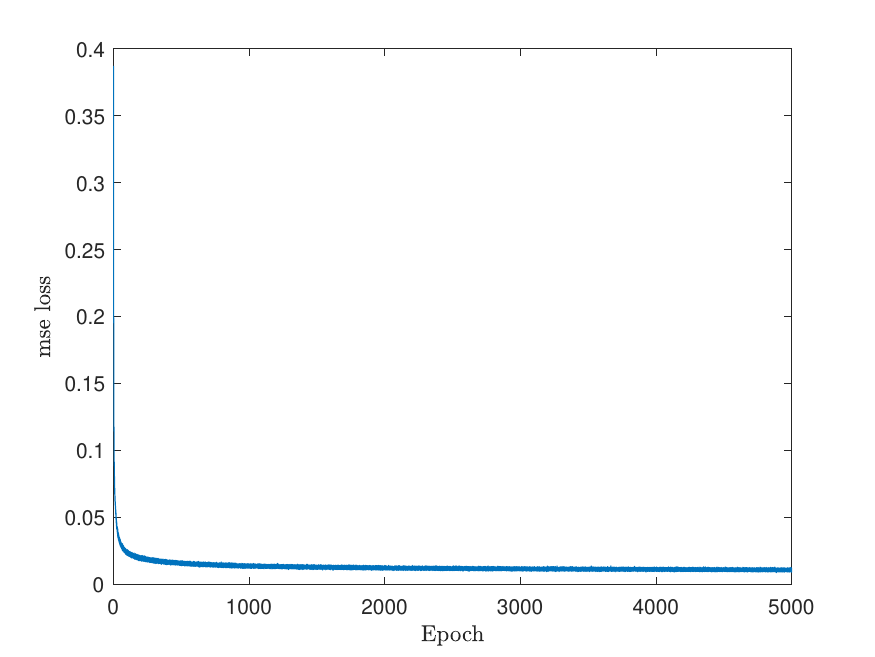}  \hspace{-0.3in}
      \includegraphics[scale=.49]{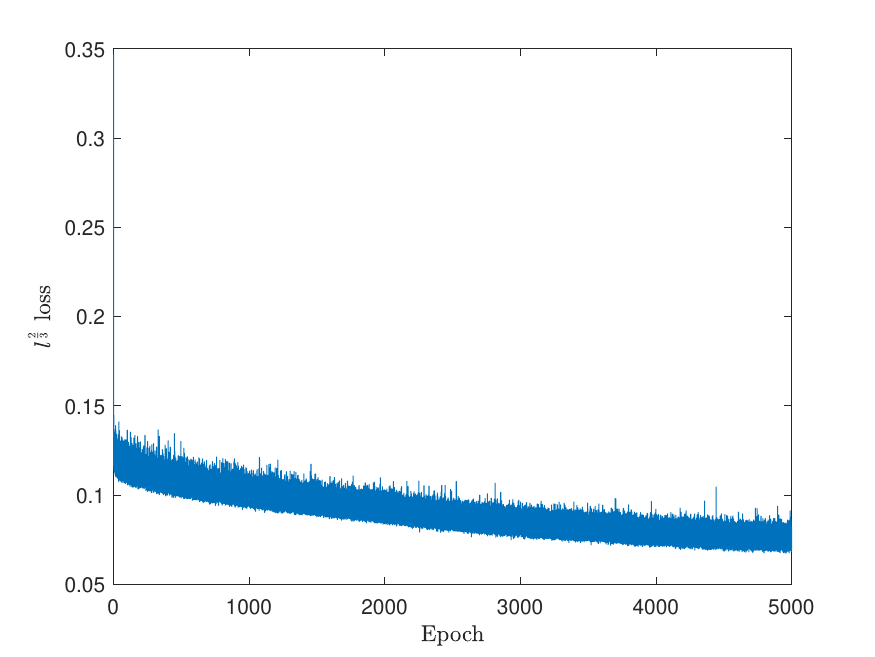}			
			       \caption{Error versus epoch number. Left: with the $l^2$ (or MSE) loss function. Right: further trained  with the
			       $l^{\f23}$ loss function. The computed network for the $l^2$ loss function was used as a starting point
			       for computing a new network with the  $l^{\f23}$ loss function.}
    \label{loss plots}
\end{figure}

In Figure \ref{loss plots}
we show an example of a run of this ADAM algorithm with 5000 epochs. 
We used a two-step training scheme; in the first step
we  used an $l^2$ (also known as MSE) loss function. In a second step, the network was further trained using 
 an $l^{\f23}$ loss function. The rationale for trying the  $l^{\f23}$ loss function is that the computed network
 seems to perform very well on average but performs poorly on some outliers 
 (Figure  \ref{error_plots}). With the  $l^{\f23}$ loss function, less importance is given to these outliers. 
 Ultimately, we only observed a modest improvement for the average reconstruction error using the network
 computed in this  two-step approach. On noisy data, the advantage of this network over the traditionally
 MSE based network disappears altogether as shown in Figure  \ref{error_plots}.



The neural network ${\cal N}_1$ was trained for geometry parameters such that 
$- \pi/2 \leq \theta \leq 0, -a_{\mbox{max}} \leq  a \leq   a_{\mbox{max}} $.
Next, we trained  a second neural
network  ${\cal N}_2$ for geometry parameters such that
$ 0 \leq \theta \leq \pi/2,
 -a_{\mbox{max}} \leq  a \leq   a_{\mbox{max}} $. 
 Separating    ${\cal N}_1$ and ${\cal N}_2$  is necessary
 due to a symmetry in the parametrization of our geometry.
 In particular, changing $m=(\theta , a$) to $m'=(\theta + \pi, -a)$
 in (\ref{bounded2 dir})--(\ref{bounded3 dir}) amounts to  changing 
 $y(t)$ to $y(-t)$
  and  $A_m \psi$ to $ - A_{m'} \tilde{\psi}$ with $\tilde{\psi} (t) = \psi(-t)$.
 Thus condition \ref{U1} fails if 
 $\theta$ has full range over $[-\pi/2, \pi/2]$.
 A third neural network ${\cal N}_3$ for $\theta$ in the range $[-\pi/2 + \epsilon, 
 \pi/2 - \epsilon] $, $\epsilon >0$  small, was trained to learn whether 
 $\theta$ is in $[-\pi/2, 0 ]$, in which case  ${\cal N}_1$  is used 
 or  $\theta$ is in $[ 0 , \pi/2]$, in which case  ${\cal N}_2$  is used.

\subsection{Testing the learned NN and stability to noise}
\label{Testing the learned neural}

We considered four types of forcing terms:
\begin{enumerate}[label=\subscript{Type \, }{{\arabic*}}]
\item The Dirichlet data $g$ in  \eqref{D2} is the incoming plane 
wave $e^{ i k x \cdot \eta}$, where $\eta$ is a unit vector in $\RR^2$.
\item The Dirichlet data $g$ in  \eqref{D2} is source point wave
$\f{i}{4} {\cal H}^1_0(k|x-s|)$ where $s \in \RR^2$ is such that 
$3 \leq |s| \leq 3.5 $. Accordingly the source is between the screen
and $S_R$.
\item The Dirichlet data $g$ in  \eqref{D2} is source point wave
  $\f{i}{4} {\cal H}^1_0(k|x-s|)$ where $s \in \RR^2$ is such that 
$5 \leq |s| \leq 7 $. Accordingly the source is outside $S_R$.
\item  This last case does not use Dirichlet problem  (\ref{D1})--(\ref{D3}).
Instead, the data $A_m \psi = u|_{S_R} $  is directly formed from 
a given forcing term.
\end{enumerate}

\begin{figure}[!ht]
   \centering
    \hspace{-0.4in}
      \includegraphics[scale=.5]{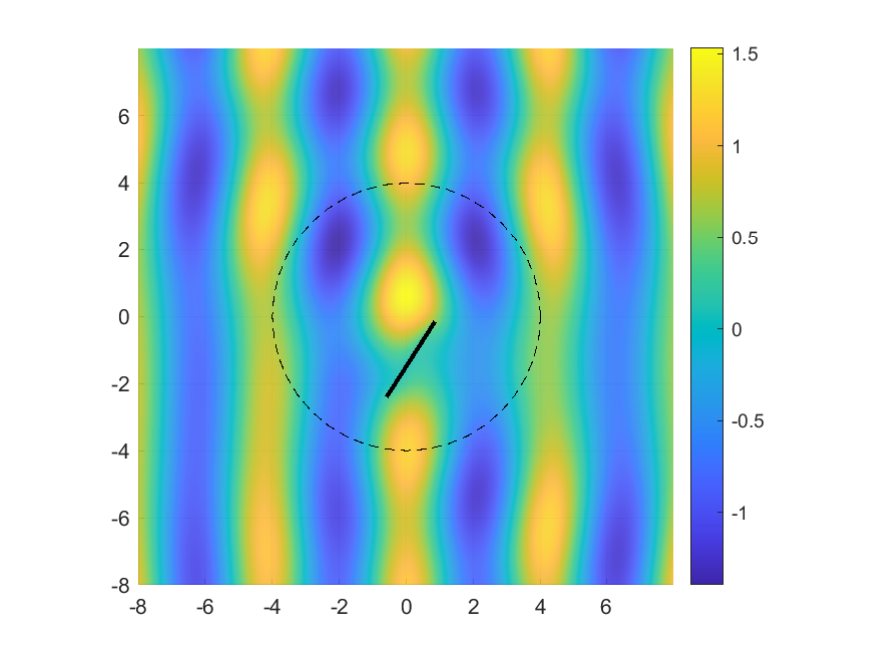} \hspace{-0.5in}
      \includegraphics[scale=.5]{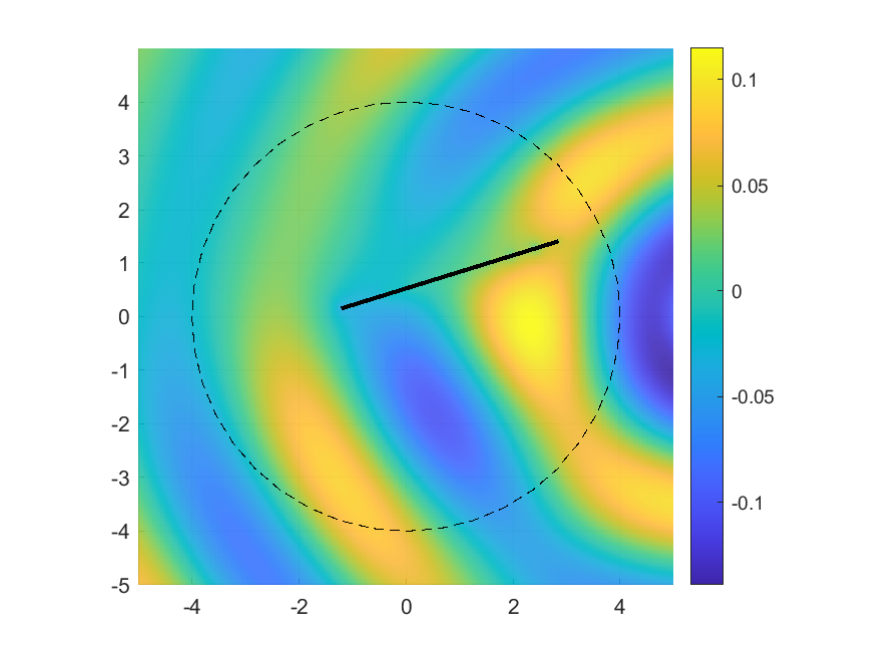}			
			       \caption{Examples 
of configurations for the fault $\Gamma$ relative to the 
circle $S_R$ in Figure \ref{fields}. The real part of the total field is visualized.
In each case, the crack $\Gamma$ is the black line segment.
The circle $S_R$ is the dotted circle. 
The left plot corresponds to a Type 1 forcing term with incidence angle $\eta = (1, 0)$. 
The right plot corresponds to Type 3 forcing term with source $s=(6,0)$. The scale is different in order 
to facilitate visualization  as
the total field quickly decays from the source. }
    \label{fields}
\end{figure}

Two examples 
 of  the fault $\Gamma$
and corresponding fields 
 are plotted in Figure \ref{fields}. 
The left plot
shows the real part of the total field
$u + e^{ i k x \cdot \eta}$ for
a Type 1 forcing term. 
The right plot corresponds to a Type 3 forcing term. The scale is different for visualization purposes as
the total field quickly decays from the source in this case.

\begin{figure}[H]
    \centering
      \includegraphics[scale=.7]{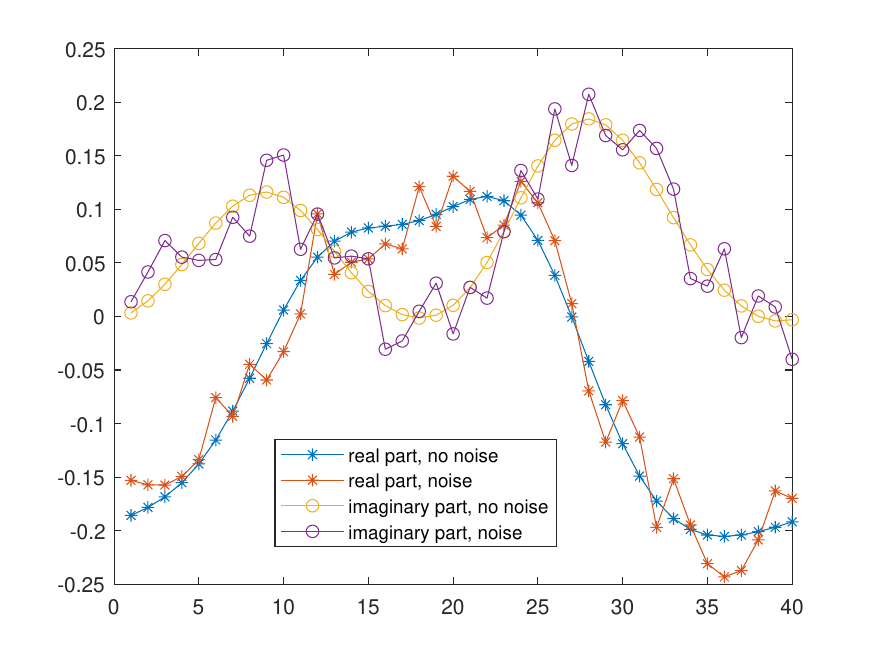}
			       \caption{Example of data in $\CC^{N_S}$.
						Real and imaginary parts of $u$ 
						at $(R \cos (j \f{2 \pi }{N_S}), 
R\sin (j \f{2 \pi }{N_S}))$, $j=1, ..., N_S $ are plotted
with $j \f{2 \pi }{N_S}$ on the horizontal axis. 
The smooth curves correspond to noise-free data. 
The jagged curves correspond to noisy data.
}
    \label{data_example}
\end{figure}

We generated data for testing the computed NN. 
The testing set comprises 1000 examples. These 1000 examples are randomly picked following
these rules (all probability distributions were taken to be uniform within their range):
\begin{enumerate}[label=\subscript{Step \, }{{\arabic*}}]
\item We randomly chose a crack boundary condition of Type 1, 2, 3, or 4 described above. 
\item If the selected type is 1, a random direction $\eta$ is picked.
 If the selected type is 2, a random source $ s \in \RR^2, 3 \leq |s| \leq 3.5$ is picked.
   If the selected type is 3, a random source $ s \in \RR^2, 5 \leq |s| \leq 7$ is picked.
    If the selected type is 4,  we pick $\psi(t) = y_1(t) -  \mathrm{i} \cos y_2(t)$ where
again $y$ depends on $t$ through equation \eqref{bounded2 dir}.
\item The geometry parameters $a, \theta$ are randomly picked within their range.
So are the parameters $o$ and $l$. This randomly defines $\Gamma$. 
\item  If the selected type is 1, 2, or 3, then the PDE (\ref{D1}-\ref{D3}) has to be solved 
to produce the data $u_{|S_R}$. 
This was done by solving integral equation \eqref{Dir int} for $[\f{\p u}{\p n}]_{| \Gamma} $
as type 1, 2, 3 only provides $u_{| \Gamma}$. In type 4, $[\f{\p u}{\p n}]_{| \Gamma} $ is directly 
given by $\psi$.
\item $u_{|S_R}$ is computed thanks to the formula $u(x) = \int_\Gamma \Phi(x,y) \left[ \f{\p u}{\p n} (y)\right] d \sigma(y)$.
It is computed at the points $(R \cos (j \f{2 \pi }{N_S}), 
R\sin (j \f{2 \pi }{N_S}))$, $j=1, ..., N_S $, $R=4$, $N_S=40$.
\end{enumerate}
To test robustness to noise, a random perturbation was added
to the data for these 1000 trials. To do that, 
for each example and each coordinate
of the corresponding data vector in $\RR^{2 N_S}$,
we drew a random number in $[-0.2, 0.2]$
which we then multiplied by 
the overall sup norm of the data and used the result 
as additive noise: see Figure \ref{data_example} for an example
of realization of clean and noisy data.

\begin{figure}[H]
   \centering
      \hspace{-0.4in} \includegraphics[scale=.49]{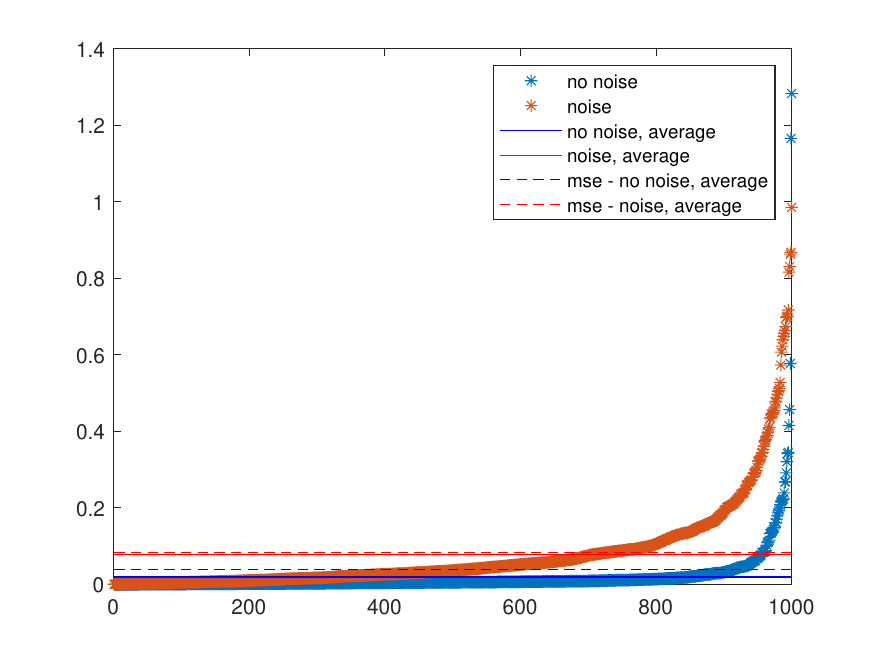}  \hspace{-0.3in}
      \includegraphics[scale=.49]{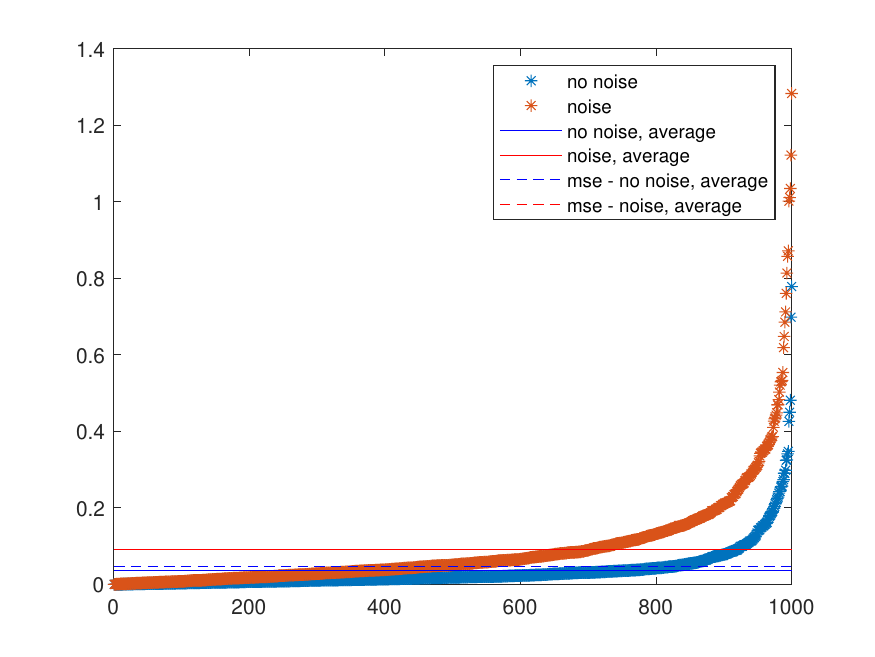}			
			       \caption{Sorted absolute value errors 
						in evaluating $\sin \theta$ (left) and $a$ (right)
						for 1000 random trials of $\theta, a$,
						support of the forcing term $g$,
						and random choice of type 1, 2, 3, or 4 forcing.
						The horizontal solid and dashed lines indicate average error for the 1000 trials using the two-step and standard (MSE) training, respectively.			
Blue: noise-free data. Red: noisy data. }
    \label{error_plots}
\end{figure}

In Figure \ref{error_plots} we plot absolute errors for computing 
$\sin \theta$ and
$a$ using networks ${\cal N}_1,  {\cal N}_2, {\cal N}_3$.
The errors are shown for 1000 randomly generated examples, with and without noise.
Collective run time for these 1000 cases is  
0.06 seconds. 
Absolute errors for $\sin \theta$ and $a$ are shown
in Figure \ref{error_plots} in blue. The average error for these 1000 
cases is  about  0.02 for $\sin \theta$ and 0.03 for $a$, for the noise-free case.
For the noisy case, these average errors rise to 0.08 and 0.09.

\subsection{Testing the learned NN with training  independent data}
   \label{finite el section}
Section \ref{Testing the learned neural}
 demonstrates the accuracy and robustness to noise of the
computed neural network for our inverse scattering problem, with synthetic data 
based on an integral representation of the scattered field,  and hence may considered to be related  to our ML algorithm 
that is based on a few singular values of an approximation to an  integral operator. 
However,  in order to avoid committing an  ``inverse crime'' 
\cite{kaipio2006statistical},
it is customary in the field of {\em inverse problems} to verify reconstruction
algorithms on   data generated    using
methods that are  entirely different from those used to build the reconstruction
algorithm. 

\begin{figure}[H]
    \centering
\hspace{-0.4in}    \includegraphics[scale=.49]{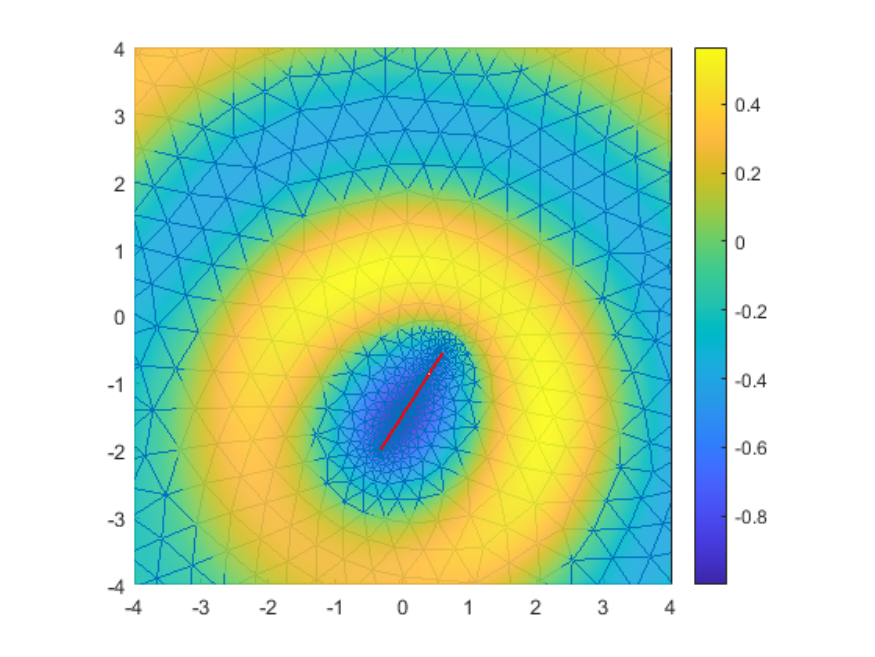} \hspace{-0.3in}			
      \includegraphics[scale=.49]{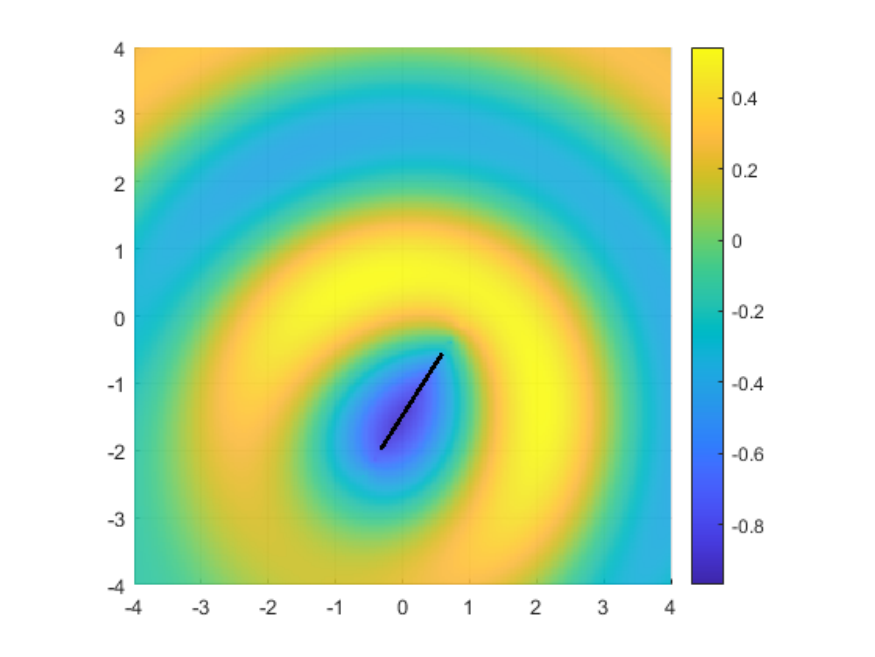}
      		       \caption{Example of  a computed scattered field $u$ using a finite element code (left) versus using integral equation
      		       \eqref{Dir int}. The scatterer is the thin rectangle (left) or the equivalent one-dimensional crack (right). Dirichlet conditions are applied in each case. The incoming wave in each case is the plane wave $e^{i k x_1}$, where $k=1.5$
      		       as in the previous cases. Note the finer mesh in the vicinity of the thin rectangle.}
    \label{finite_el}
%
\end{figure}

To this end, in this section for the learned NN testing experiments,   we  generate synthetic scattering data $u$ on $S_R$
using a finite element method (FEM) solver that directly solves the PDE in a bounded domain by approximating the radiation condition~\eqref{D3} and
modeling the crack $\Gamma$ as a thin rectangle.
The decay condition  \eqref{D3} is replaced by the absorbing condition 
$ \f{\p u}{\p r} - i k u = 0 $ on a large circle with radius $R_2 = 20$.
The related FEM-based  data is thus less accurate and independent than the previously used 
integral equation based solver: this is also  partly due to the need of using a very fine
mesh close to the thin rectangle as shown in~Figure  \ref{finite_el}.
However, FEM solver based approach presents the advantage 
of providing a more realistic model where domains are not infinite 
and cracks are not infinitely thin. 

We then produce 1000 randomly generated Type 1 cases (plane wave impinging on a crack) 
with the parameters $a, \theta, o, l$ sampled with their range (using uniform distributions) as the incidence angle $\eta $.
Each case is done using the FEM: the cracks are truly very thin rectangles
(with thickness 0.01).
As previously, the scattered field $u$ is sampled at the points
$(R \cos (j \f{2 \pi }{N_S}), 
R\sin (j \f{2 \pi }{N_S}))$, $j=1, ..., N_S $, $R=4$, $N_S=40$.
That being said these points are not necessarily vertices of the mesh used by the FEM solver.  The value of the solution
at those points is inferred by interpolation. This is yet another difference between this FEM  and 
the integral representation for data generation in the previous section.

\begin{figure}[H]
    \centering
 \hspace{-0.45in}     \includegraphics[scale=.49]{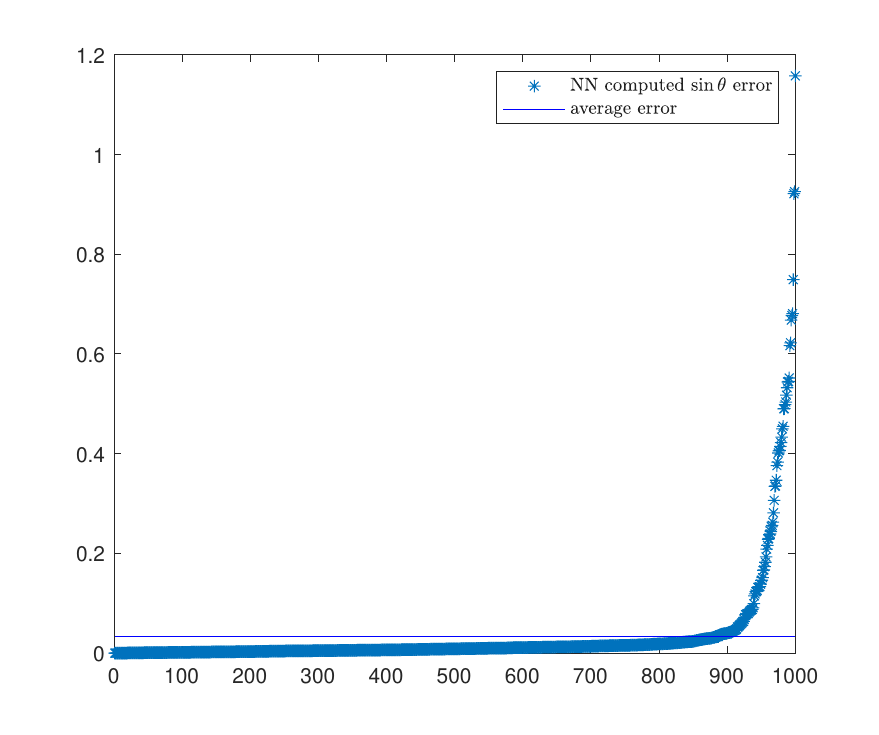} \hspace{-0.28in}
       \includegraphics[scale=.49]{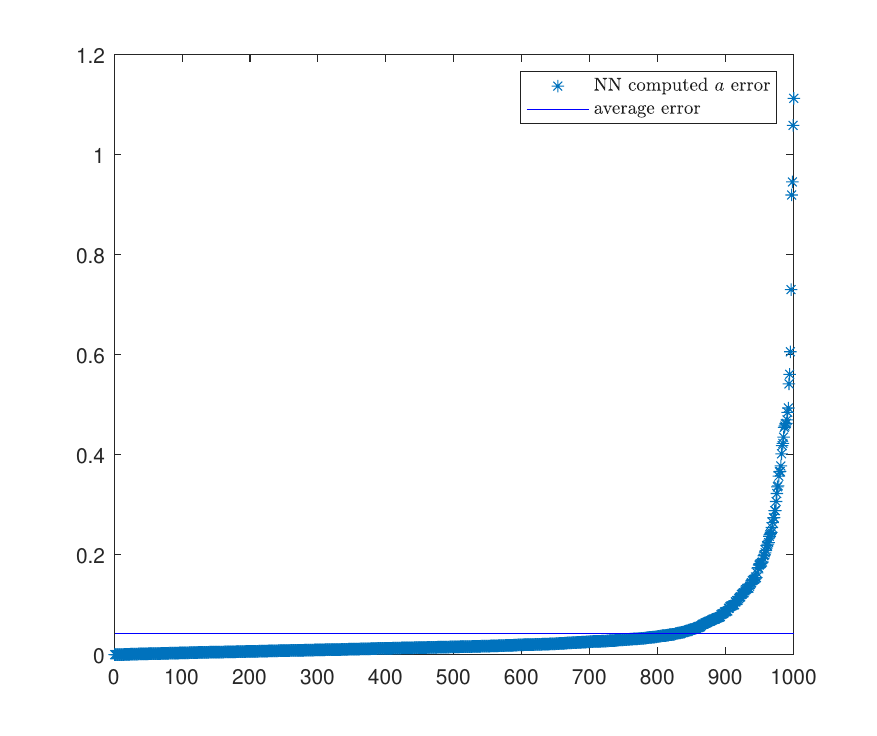}			
			       \caption{Sorted absolute value errors 
						in evaluating $\sin \theta$ (left) and $a$ (right)
						for 1000 random trials of $\theta, a, o,  l, \eta$,
						in the case of an incoming plane wave $e^{i k x \eta}$
						with Dirichlet conditions on the crack.
						Here, the data for the inverse problem was generated by a finite element solver.
						The neural networks ${\cal N}_1, {\cal N}_2, {\cal N}_3$ were used to evaluate 
						$\sin \theta$ and $a$.
						The horizontal solid lines indicate average error for the 1000 trials.
}
%
    \label{NN_with_finite_el}
\end{figure}

 The sampled data is then fed in the neural networks  ${\cal N}_1, {\cal N}_2, {\cal N}_3$  to evaluate 
						$\sin \theta$ and $a$. Results are shown in Figure \ref{NN_with_finite_el}.
 The mean error was found to be 0.03 for $\sin \theta$ and 0.04
for $a$. 
This is very close to the noise-free errors found in section \ref{Testing the learned neural} and better than the errors found in
the noisy case.

\section{Appendix: proofs of stated results in  sections \ref{prelim}  \& \ref{Dirichlet screen}}
We first  recall Theorem~2.2 from~\cite{alberti2023inverse}, incorporating the manifold and differentiability concepts introduced in Definitions~2.1--2.3 of \cite{alberti2023inverse} adapted to our setting.
The finite-dimensional differentiable manifold $\mathcal{M}$
will be in our case a subset of ${\cal B}' \times ( E\setminus \{ 0\}) $
which is itself in the Banach space $\RR^p  \times E$.
\begin{thm} \label{Alberti thm}
Let $\mathcal{M}$ be a finite-dimensional manifold  of ${\cal B}' \times (E\setminus \{ 0 \})$.
Assume that ${\bf \Psi}$ is injective on $\mathcal{M}$ and that 
the differential of ${\bf \Psi}$ is injective at every point in $\mathcal{M}$.
Let ${\rm \bf K}$ be a compact subset of $\mathcal{M}$.
Then there is a positive constant $C$ such that for all
$(m, \varphi), (m',\psi) $ in  ${\rm \bf K}$,
\bean \label{Alberti form}
 \| {\bf \Psi}(m,\varphi) - {\bf \Psi}(m',\psi) \| 
 \geq C  (| m - m'| + \| \varphi - \psi\| ).
\eean
\end{thm}
Our proof of Lipschitz stability involves constructing a manifold composed of pairs $(m, \varphi)$, where:
\begin{itemize}
\item $m$ lies in an open subset of parameters;
\item $\varphi$ belongs to the range of a spectral projector $P_m$ associated with the  operator $A_m^* A_m$;
\item the mapping $m \mapsto P_m$ is continuous within this open set.
\end{itemize}
A direct application of Theorem~\ref{Alberti thm} yields local Lipschitz continuity of the inverse map ${\bf \Psi}$. To extend this continuity globally over ${\cal B}$, a more delicate argument is required, as the spectral projectors $P_m$ may not vary continuously with $m$ outside the local neighborhood \cite{kato2013perturbation}.
The finite-dimensional subspaces defined by the ranges of the $P_m$ projectors serve a dual purpose: they support the mathematical analysis and also facilitate the learnability of the manifold structure in our ML-motivated framework.

\vspace{0.1in}
\noindent
\textbf{Proof of Proposition \ref{U1 U2 eq}}\\
Condition \ref{U1} is clearly equivalent to the injectivity of ${\bf \Psi }$ on ${\cal B}'\times (E \setminus \{ 0 \})$.\\
The derivative of ${\bf \Psi }$ at $(m,\varphi) \in {\cal B}' \times E $ is a continuous linear function from
$\RR^p \times E$ to $F$ which can be simplified due to the particular form   \eqref{PSsi def}. 
 Let $e_1, ..., e_p$ be the natural basis of $\RR^p$, $a_1, ... , a_p \in \RR$ and $\psi \in E$. We note that
 \bea
 &&{\bf \Psi} (m + \sum_{j=1}^p a_j e_j , \varphi + \psi) - {\bf \Psi}(m,\varphi) \\
& =& {\bf \Psi} (m + \sum_{j=1}^p a_j e_j , \varphi + \psi) - {\bf \Psi}(m,\varphi+\psi) 
+  {\bf \Psi}(m,\varphi+\psi)  - {\bf \Psi}(m,\varphi) \\
 &=&  \sum_{j=1}^p  a_j \p_{e_j} A_m (\varphi+\psi) + o(\sum_{j=1}^p |a_j| )+ A_m \psi\\
&=&  \sum_{j=1}^p  a_j \p_{e_j} A_m \varphi + \sum_{j=1}^p  a_j \p_{e_j} A_m \psi
+ A_m \psi+ o(\sum_{j=1}^p |a_j|) .
 \eea
Since 
\bea
\| \sum_{j=1}^p  a_j \p_{e_j} A_m \psi \| \leq \f12  
\sum_{j=1}^p | a_j |^2 + \| \p_{e_j}  A_m \psi \|^2,
\eea
this term is $o(\sum_{j=1}^p |a_j| + \| \psi\|)$. Altogether, we have found that the derivative of 
${\bf \Psi }$ at $(m,\varphi)$ is the linear function
\bea
\RR^p \times E \ri F \\
(\sum_{j=1}^p a_j e_j, \psi) \mapsto  \sum_{j=1}^p  a_j \p_{e_j} A_m \varphi + A_m \psi .
\eea
Assume now that 
the derivative of ${\rm \bf \Psi }$ is injective at all points in 
${\cal B}' \times (E \setminus \{ 0\})$.
Let $q=(q_1, ..., q_p) \in \RR^p$ be such that $|q| =1$ and $\varphi \neq 0 $ in $E$. Let $\psi$ be in $E$. Then 
 $ \sum_{j=1}^p  q_j \p_{e_j} A_m \varphi + A_m \psi \neq 0$, since $q \neq 0$. 
 Now assume that $\varphi =0,  \psi\neq 0$. Then $ \sum_{j=1}^p  q_j \p_{e_j} A_m \varphi + A_m \psi  = A_m \psi \neq 0$, 
 as ${\rm \bf \Psi }$ is injective on
 ${\cal B}' \times (E \setminus \{ 0\})$.
 We have shown that \ref{U2} holds.\\
 Conversely, assume that \ref{U2} and \ref{U1} hold. Then if $\varphi\neq 0$ and $\sum_{j=1}^p a_j \neq 0$ or $\psi \neq 0 $,
 then   $\sum_{j=1}^p  a_j \p_{e_j} A_m \varphi + A_m \psi  \neq 0$. Thus the  derivative of 
 ${\bf \Psi }$ at $(m,\varphi) \in {\cal B}' \times (E\setminus \{ 0 \})$ is injective.
 $\square$ \\\\
\textbf{Proof of Theorem \ref{main theorem}}\\
We start by introducing orthogonal projections 
	on the first few eigenspaces of $A_m^* A_m$.
	For  $m \in {\cal B}'$ let 
	\bea
	 \lambda_1^2(m) > ... > \lambda_N^2 (m)>  ...
	\eea
	be the distinct ordered eigenvalues of $A_m^*A_m$. 
	Fix $m_1 \in {\cal B}'$ and
	let ${\cal C}_1$ be the  circle in the complex plane centered at the origin with radius
 $\max_{m \in {\cal B}} \| A_{m}^* A_{m}\|+1$, and ${\cal C}_2$ be the  circle  centered at the origin with radius
 $\ds \f{\lambda_N^2(m_1)+\lambda_{N+1}^2(m_1)}{2}$.
 
For $m$ near $m_1$ define  \cite{kato2013perturbation},
\bean \label{int for proj}
P_{m,1} = \f{1}{2\mathrm{i} \pi} \int_{{\cal C}_1} (zI -  A_{m}^* A_{m} )^{-1} \, dz -
 \f{1}{2\mathrm{i}\pi} \int_{{\cal C}_2} (zI -  A_{m}^* A_{m} )^{-1}  \, dz . 
\eean
Note that for $m$ near $m_1$,  $P_{m,1}$ is also the 
orthogonal projection 
on the sum of the eigenspaces
	of $A_{m}^*A_{m}$ corresponding to the eigenvalues 
	greater than $\ds \f{\lambda_N^2(m_1)+\lambda_{N+1}^2(m_1)}{2}$.
	Based on formula  \eqref{int for proj},
	we can argue that 
		$m \mapsto P_{m,1}$ is a  $C^1$ function for 
	$m$ in an open ball $B(m_1, 3 \epsilon)$
        centred at $m_1$
        with radius $3 \epsilon>0$.
	Necessarily, $s_1=\dim R(P_{m_1,1})$ is constant in $B(m_1, 3 \epsilon)$ by 
	continuity of the trace.
	
	We claim that, after possibly shrinking $\epsilon$,
	there is a linear bijection $b_{m,1}: \RR^{s_1} \to R(P_{m,1})$
	such that $m \mapsto b_{m,1}$ is a $C^1$ function  in $B(m_1, 3 \epsilon)$.
	Indeed, let $\varphi_{1,1}, ..., \varphi_{1,s_1}$ be an orthonormal basis of $R(P_{m_1,1})$.
	The matrix $\langle P_{m_1,1} \varphi_{1,i}, P_{m_1,1} \varphi_{1,j} \rangle_{1 \leq i,j \leq s_1}$ is the identity matrix since
	$P_{m_1,1} \varphi_{1,i}=\varphi_{1,i}$. By continuity,  
	$\langle P_{m,1} \varphi_{1,i}, P_{m,1} \varphi_{1,j} \rangle_{1 \leq i,j \leq s_1}$
	is invertible for $m$ near $m_1$, proving that 
	$P_{m,1}\varphi_{1,1}, ..., P_{m,1}\varphi_{1,s_1}$
	is a basis of $ R(P_{m,1})$ since its dimension is also $s_1$.
	
	Next, we cover ${\cal B}$ by finitely many  balls
	 $B(m_i, \epsilon)$, $i \in I$.  Choosing $\epsilon >0$
	 small enough,  
	 we can assume that $B(m_i,  3 \epsilon)\subset {\cal B}' \subset \RR^p,~ i \in I$,
	and  that $P_{m,i}\varphi_{i,1}, ..., P_{m,i}\varphi_{i,s_i}$
	is a basis of $ R(P_{m,i})$, if $m$ is in $B(m_i,  3 \epsilon)$.
	Finally,  an explicit formula for $b_{m,i}$ with $m$ in $B(m_i,  3 \epsilon)$
	can be given.
	For $m \in B(m_i,  3 \epsilon)$, $i \in I$ we set
	\bea
	b_{m,i}: \RR^{s_i}  \ri R(P_{m,i}), \\
	(a_1, ..., a_{s_i}) \ri   a_1 P_{m,i} \varphi_{1,i}+  ... + a_{s_i}  P_{m,i}  \varphi_{s_i,i}.
	\eea
	That way we have defined a singular vectors-based manifold 
	\bea \mathcal{M}_i = \{ (m, \varphi): m \in B(m_i, 3 \epsilon ) , \varphi \in R(P_{m,i}) \}
	\eea
	with dimension $ p + s_i$, for $i \in I$. 
	\begin{lem}\label{proj case}
Assume that conditions \ref{U1} and \ref{U2} hold.
There is a positive constant $C$ such that 
for all $m,m' \in {\cal B}$,
\bean \label{stability proj case}
\| A_m \varphi - A_{m'}  \psi  \| \geq 
C(\| \psi\| |m-m'| + \| \varphi -\psi\|),
\eean
for all $\varphi\in R(P_{m,i})$,
$\psi \in R(P_{m,j})$ 
where $i, j \in I $ are such that 
 $m \in B(m_i, \epsilon) $, $m' \in B(m_j, \epsilon )$. 
\end{lem}

\vspace{0.01in}
\noindent
\textbf{Proof of Lemma \ref{proj case}}\\
Fix $i \in I$.
We can apply Theorem \ref{Alberti thm} 
in $M $
and the compact subset  of $M$,
${\rm \bf K} = \{  (m,\varphi) \in M: m \in \ov{B({m_i}, 2 \epsilon)}
, \varphi \in  R(P_{m,i}) , \| \varphi \| =1\}$.
As a result, there is a constant $C_i>0$ such that
	\bean
\| A_m \varphi- A_{m'} \psi  \| \geq 
C_i(|m-m'| + \|  \varphi - \psi \|),
\label{Wmi}
\eean
for all $m, m' \in \ov{B({m_i}, 2 \epsilon)}$, $\varphi \in R(P_{m,i})$,
$\psi\in R(P_{m',i})$ with  $\| \varphi \| = \| \psi \| =1$.
By homogeneity,  
\bea
\| A_m \|\psi \|\varphi- A_{m'} \psi  \| \geq 
C_i(\|\psi \| |m-m'| + \|  \|\psi\|\varphi - \psi \|),
\eea
for any $\psi  \in R(P_{m',i})$ and $\varphi \in R(P_{m,i})$,  with  $\| \varphi \|=1$,
so using that $ R(P_{m,i})$ is a linear space,
\bea
\| A_m \varphi- A_{m'} \psi  \| \geq 
C_i(\|\psi \| |m-m'| + \|  \varphi - \psi \|),
\eea
for any $\psi  \in R(P_{m',i})$ and $\varphi \in R(P_{m,i})$, and any $m, m' \in 
\ov{B({m_i}, 2 \epsilon)}$.
Setting $C=\min_{i \in I }C_i$, 
it follows that
\bea
\| A_m \varphi- A_{m'} \psi  \| \geq 
C (\| \psi\| |m-m'|  + \|  \varphi - \psi\|),
\eea
for all $\psi  \in R(P_{m',i})$, $\varphi \in R(P_{m,i})$, $i \in I$, and $m, m' \in {\cal B}$ such that $ |m - m_i| \leq 2\epsilon$,
$ |m' - m_i| \leq 2\epsilon$.
Next assume, that $m,m'$ are such that $|m-m'| \leq  \epsilon$, $m \in 
B(m_i, \epsilon)$,
$m' \in B(m_j, \epsilon)$, but $i \neq j$.
Let $\varphi \in R(P_{m,i})$, $\psi \in R(P_{m',j})$.
Observe that 
$R(P_{m,i}) \subset R(P_{m,j})$
or $R(P_{m,j}) \subset R(P_{m,i})$. Assume
that $R(P_{m,i}) \subset R(P_{m,j})$.
Then by \eqref{Wmi},  as $m$ and $m'$ are in $B(m_j, 2 \epsilon)$,
\bea
\| A_m P_{m,j} \varphi- A_{m'} P_{m',j} \psi  \| \geq 
C(\|P_{m',j} \psi \|  |m-m'| + \|  P_{m,j} \varphi - P_{m',j} \psi \|).
\eea
But as $R(P_{m,i}) \subset R(P_{m,j})$, $P_{m,j} \varphi= \varphi$. 
Since $P_{m',j} \psi  =\psi$, 
 we find that
\bea
\| A_m \varphi- A_{m'} \psi  \| \geq 
C(\|\psi \|  |m-m'| + \|  \varphi - \psi \|).
\eea
The case $R(P_{m,j}) \subset R(P_{m,i})$ is handled similarly. $\square$\\
To finish proving Theorem \ref{Alberti thm}, 
it remains  to show that 
\bean
\inf
\f{\| A_m  P_{m,i}  \varphi - A_{m'} P_{m',j}\psi \|}{\| P_{m',j}\psi\| |m-m'| + 
\|  P_{m,i} \varphi -P_{m',v}\psi\|} >0, \label{big inf}
\eean
where the infimum is taken over the set $i,j \in I$, $m \in B(m_i, \epsilon), m' 
\in B(m_j, \epsilon)$,
such that $|m - m'|\geq \epsilon$ and $\varphi, \psi \in E$ such that the denominator
$ \| P_{m',j}\psi\| |m-m'|+ 
\|  P_{m,i} \varphi -P_{m',v}\psi\|$ is positive.
If $P_{m',j}\psi=0$ it suffices to show that 
\bean \label{suffices}
\inf_{  i \in I, \varphi \in B(m_i, \epsilon), \|  P_{m,i} \varphi\|>0}
\f{\| A_m  P_{m ,i} \varphi\|}{ \|  P_{m,i} \varphi \|} >0.
\eean
By \eqref{int for proj}, $P_m$ is the orthogonal projection on the sum of eigenspaces
corresponding to the  eigenvalues of 
$A_m^*A_m$
greater than $\f{\lambda_N^2(m_i)+\lambda_{N+1}^2(m_i)}{2}$
so the infimum in \eqref{suffices} is greater than \\
$\min_{i \in I} \left[\f{\lambda_N^2(m_i)+\lambda_{N+1}^2(m_i)}{2}\right]^{\f12}$.\\
Now, let $d$ be greater than the diameter of ${\cal B}$. Since 
 $\|P_{m',j} \psi\| |m-m'| + \| P_{m,i} \varphi -P_{m',j}\psi\| \leq (d+1)\| P_{m',j}\psi\| 
 + \| P_{m,i}  \varphi \| $, it suffices to show that
 \bea
\inf
\f{\| A_m P_{m,i} \varphi - A_{m'}P_{m',j} \psi \|}{(d+1) \| P_{m',j}\psi\| + \|P_{m,i} \varphi \|} > 0,
\eea 
where the infimum is taken over the same set as in the infimum in \eqref{big inf}.
Arguing by contradiction, assume that there are two sequences 
$(i_n), (j_n) \in I$, 
$(m_n), (m_n') \in { \cal B}$ with $m_n \in B(m_{i_n}, \epsilon), m_n' \in B(m_{j_n},
\epsilon)$ such that $|m_n - m_n'|\geq \epsilon$, and two sequences
$(\varphi_n), (\psi_n) \in E$  such that 
$\|P_{m_n', j_n} v_n\|>0$ and
\bean \label{lim formula}
\lim_{n \ri \infty}
\f{\| A_{m_n} P_{m_n, i_n} \varphi_n - A_{m_n'} P_{m_n', j_n} \psi_n \|}{(d+1)
 \| P_{m_n', j_n}\psi_n\| 
+ \| P_{m_n, i_n}\varphi_n \|} = 0.
\eean
Set $\omega_n = \f{ \| P_{m_n, i_n}\varphi_n\| }{(d+1) 
\| P_{m_n', j_n} \psi_n\| + \|P_{m_n, i_n} \varphi_n\|}$.
We have from \eqref{lim formula}
\bea
\lim_{n \ri \infty} \omega_n  A_{m_n}  \f{ P_{m_n, i_n} \varphi_n}{\| P_{m_n, i_n} \varphi_n \|}
- (1- \omega_n) (d+1)^{-1}A_{m_n', j_n}\f{  P_{m_n', j_n}\psi_n}{\| P_{m_n', j_n}\psi_n \|}  =0.
\eea
By compactness, after extracting subsequences
we can assume that  
$i_n \ri i, j_n \ri j$ in $I$, 
$\omega_n $ converges to $\omega \in [0,1]$,
$m_n$ converges to $m \in \ov{B(m_i, \epsilon )}$, and $m_n'$ converges to $m' \in \ov{B(m_j, \epsilon)}$. 
Necessarily  $|m - m'| \geq \epsilon$. Note that  
$\f{ P_{m_n, i } \varphi_n}{\| P_{m_n ,i } \varphi_n \|} =P_{m_n, i } \f{ P_{m_n, i} \varphi_n}{\| P_{m_n, i} \varphi_n \|}$.
Since each operator $P_{m,i}$ is compact and $m \mapsto P_{m,i}$ is continuous in the closed ball
$\ov{B(m_i,  \epsilon)}$, after extracting a subsequence,
$\f{ P_{m_n, i } \varphi_n}{\| P_{m_n ,i } \varphi_n \|}$ converges strongly to some $\tilde{\varphi} \in R(P_{m,i})$ with $\|\tilde{\varphi}\|=1$. 
Similarly,  $\f{  P_{m_n'}\psi_n}{\| P_{m_n'}\psi_n \|} $ converges strongly to some
$\tilde{\psi} \in R(P_{m',j})$ with $\|\tilde{\psi}\|=1$.
At the limit we find
\bea
\omega  A_{m}  \tilde{\varphi}
- (1- \omega) (d+1)^{-1}A_{m'} \tilde{\psi}  =0,
\eea
contradicting condition  \ref{U1}.  $\square$ \\\\

Thanks to Lemma \ref{proj case}, we can now finish proving 
Theorem \ref{main theorem}.
Constructing $B(m_i, \epsilon), i \in I$ as in the proof of 
Lemma \ref{proj case}, $s_i $, which denoted the dimension of 
$R(P_{m,i})$ satisfies $s_i \geq N$. It follows that 
$E_{m,N} \subset R(P_{m,i}) $, for all $i \in I$ and $m$ in $B(m_i, \epsilon)$.
Since the balls $B(m_i, \epsilon)$ cover ${\cal B}$, the result is proved. $\square$\\\\


\vspace{0.1in}
\noindent
\textbf{Proof of Proposition \ref{direct  dir prob}}\\
We first show uniqueness. Assume that $g=0$.
Let $S_R$ be the sphere centered at the origin with radius $R$. 
Applying Green's theorem
we find that $\mbox{Im} \int_{S_R} u \f{\p \ov{u}}{\p r} = 0$. 
Next, since $u \in {\cal V}$, 
there is a sequence $R_n \ri \infty$ such that,
\bea
\lim_{n \ri \infty}\int_{S_{R_n}} \left|\f{\p u}{\p r} - i k_0 u\right|^2 = 0,
\eea 
so altogether we have
\bea
\lim_{n \ri \infty}\int_{S_{R_n}} \left|\f{\p u}{\p r}\right|^2 + | u|^2 = 0.
\eea 
Due to Rellich's lemma for far field patterns, it follows that $u (x) =0$, if $|x| >R_0$. 
Since the only regularity assumption on  $k$ is that it is in $L^\infty$, 
there is no elementary argument for showing that $u(x)$ is zero if $|x| \leq R_0$.
However, we can use results from the unique continuation literature,
in particular, the corollary of  Theorem 1 in \cite{barcelo1988weighted}, 
to claim that $u$ is zero throughout $\RR^d\setminus \ov{\Gamma}$.\\
Next, we show existence. 
Fix $R' >R_0$ and let $B_{R'}\subseteq \RR^d$ denote the open ball
with radius $R'$ centered at the origin.
We can extend $g$ to a function 
$\psi
\in H^1(\RR^d)$
supported strictly inside $B_{R'}$. 
Using a continuous extension operator,  $\psi$ depends continuously 
on $g$.
We now seek to solve an equivalent problem for $\tilde{u} =u - \psi$.
Define the closed subspace 
$H^1_{\Gamma, 0}(B_{R'}) \subseteq H^1(B_{R'}) $,
\bea
H^1_{\Gamma, 0}(B_{R'}) = \{ w \in H^1(B_{R'}): v=0 \mbox{ on } \Gamma \}.
\eea
Note that this definition requires the trace of $v$ on $\Gamma$ to be zero
on each side.
 Define the bilinear functional,
\bean \label{B def}
{\mathrm b}(v,w) = \int_{B_{R'}} \nabla v \cdot \nabla w -  k^2  v w  - \int_{S_{R'}} T_{R', k_0} v w,
\eean
for $ v, w \in H^1(B_{R'}) $ and where $T_{R',k_0}$ is the Dirichlet to Neumann map for radiating solutions to the Helmholtz equation in the exterior of $B_{R'}$
with wavenumber $k_0$. $T_{R',k_0}$ is known to be a  continuous
mapping from $H^{\f12 } (S_{R'})$ to $H^{-\f12 } (S_{R'})$, 
while $-T_{R',0}$ is strictly coercive, and $T_{R',k_0} - T_{R',0}$ is compact
from $H^{\f12 } (S_{R'})$ to $H^{-\f12 } (S_{R'})$, see 
Section~5.2 of~\cite{colton1998inverse} or Section~2.6.5 of~\cite{nedelec2001acoustic}. 
According to the uniqueness property covered above, we have that if
$v \in H^1(B_{R'})$ and
${\mathrm B}(v,w) =0 $ for all $w \in H^1(B_{R'})$, then $v=0$. \\ 
Now, consider the variational problem:
\bean \label{var pb dir}
\mbox{find } \tilde{u} \in H^1_{\Gamma, 0}(B_{R'}) \mbox{ such that } \forall w \in 
H^1_{\Gamma, 0}(B_{R'}), \no \\
{\mathrm b}(\tilde{u},w) = -{\mathrm b}(\psi,w).
\eean
This problem has at most 
one solution since ${\mathrm b}$ is non-degenerate. 
Existence follows by arguing that this problem is in the form strictly coercive plus compact,
which follows from the properties of the operator $T_{R',k_0}$ recalled above.
Finally, $u = \tilde{u} + \psi$ 
can be extended to $\RR^d \setminus \ov{\Gamma}$
as a function satisfying 
  (\ref{D1}-\ref{D3}).
$\square$

\vspace{0.1in}
\noindent
\textbf{Proof of Theorem \ref{InverseProblemResultDir}}\\
Set 
$U=\RR^{d} \setminus \ov{\Gamma_1 \cup \Gamma_2}$, $u= u^1 - u^2 \in U$.  We can then argue as in the proof of Proposition 
\ref{direct  dir prob}
that $u$ is zero in $U$.
Next,  we argue by contradiction: suppose that there is an $x \in \Gamma_1 \setminus  \ov{\Gamma_2}$. Then there is an open ball $B(x,r)$ centered at $x$
with radius $r>0$ such that $B(x,r) \cap \ov{\Gamma_2} = \emptyset$.
Now $\left[ \f{\p u}{\p n} \right]=0$ on $B(x,r) \cap \ov{\Gamma_1} $,
and as $(\Delta +k ^2)u^2 = 0$ in $B(x,r)$,
$\left[ \f{\p u^2}{\p n} \right]=0$ on $B(x,r) \cap \ov{\Gamma_1} $.
It follows that $\left[ \f{\p u^1}{\p n} \right]=0$ on $B(x,r) \cap \ov{\Gamma_1} $,
contradicting the assumption  that $\ds \left[ \f{\p u}{\p n} \right]$ has full support in $\Gamma$. We conclude that 
$\Gamma_1 \subset \ov{\Gamma_2}$. Reversing the roles of $\Gamma_1$ and 
$\Gamma_2$ we then find that $\ov{\Gamma_1} =\ov{\Gamma_2}$. Using one more time
that $u$ is zero in $U$,
since $u=0$ on each side of  $\Gamma_1= \Gamma_2$, it follows that
$g_1- g_2 =0$ almost everywhere in $\Gamma_1$.
$\square$

\vspace{0.1in}
\noindent
\textbf{Proof  of Proposition \ref{BU2 cond}}\\
To prove the first point, 
we set $u^1(x)=\int_{-M}^{M} \Phi(x,y) \psi (t) dt$,
$u^2(x)=\int_{-M}^{M} \Phi(x,y) \phi (t) dt$, with $y$ 
as in \eqref{bounded2 dir}.
Then by Theorem \ref{InverseProblemResultDir}, 
$m=m'$ and $u^1=u^2$ outside the line defined by \eqref{bounded2 dir}.
It follows that the jump of the normal derivative of $u^1$ across that line 
equals the jump of the normal derivative of $u^1$ so $\psi=\phi$.\\
To prove the second point, 
assume that $\p_{q} A_m \psi - A_m \phi=0 $ for some 
$\psi, \phi  \in H^{-\f12 }((-M,M))$, and $q=(q_1,q_2)$,  $|q|=1$.
According to the chain rule, 
for all $x \in S_R$, 
\bea
&&\p_{q} A_m \psi -A_m \phi \\
&=& q_1 \f{\p}{\p \theta  } A_m \psi + 
q_2 \f{\p}{\p a  } A_m \psi -A_m \phi \\
 &=&   q_1 \int_{-M}^{M}  \nabla_y \Phi(x,y) 
\cdot (nt - a \tau)  \psi(t)  dt
+ q_2 \int_{-M}^{M}  \nabla_y \Phi(x,y) 
\cdot n \psi(t)  dt  \\
&& -  \int_{-M}^{M} \Phi(x,y) 
\phi(t)  dt .
\eea
Define  $w(x)$   by the formula in the previous line  
 for all
$x \in \RR^2 \setminus \ov{\Gamma}$.
By construction, $w(x)=0$
for all $x \in S_R$, and $(\Delta + k^2)w=0$
in $\RR^2 \setminus \ov{\Gamma}$.
Since 
$\f{w}{\sqrt{1+r^2}\ln (2 + r)}, 
\f{\nabla w}{\sqrt{1+r^2} \ln (2 + r)},
 \f{\p w}{\p r} - i k w\in L^2(\RR^{2}\setminus \ov{\Gamma})$,
 it follows that $w$ is zero 
in $\RR^2 \setminus \ov{\Gamma}$.
Now, the term $ \int_{-M}^{M} \Phi(x,y) 
\phi(t)  dt $ is known to be continuous  across $\Gamma$. 
The jump of $w$  across $\Gamma$ can be determined according to the rules
shown in Lemma~1 or Appendix~C of~\cite{volkov2021stability}. 
According to these jump formulas 
we find that 
$(q_1 t  + q_2)\psi(t)  =0 $
inside the support of $\psi $, so $\psi=0$.
At this stage, we just need to notice that since $\psi=0$, the jump 
of the normal derivative  of $w$ across $\Gamma$ is $\phi$, so $\phi=0$.
$\square$

\textbf{Acknowledgements}

SCH and MG gratefully acknowledge the support of
the Australian Research Council (ARC) Discovery Project Grant (DP220102243).
MG is supported by the Simons Foundation through Grant No. 518882.
 DV is supported by the Simons Foundation through Grant MPS-TSM-00007534.



\end{document}